\documentclass[10pt,twoside,a4paper,reqno]{amsart}
\usepackage{amscd,amsmath,amsthm,amsfonts,latexsym,amssymb}

\theoremstyle{plain}
\newtheorem{theo+}           {Theorem}
\newtheorem{prop+}           {Proposition}
\newtheorem{coro+}           {Corollary}
\newtheorem{lemm+}           {Lemma}
\newtheorem{conjecture}   {Conjecture}

\theoremstyle{definition}
\newtheorem{rema+}           {Remark}
\newtheorem{defi+}           {Definition}

\newenvironment{theorem}{\begin{theo+}}{\end{theo+}}
\newenvironment{proposition}{\begin{prop+}}{\end{prop+}}
\newenvironment{corollary}{\begin{coro+}}{\end{coro+}}

\newenvironment{remark}{\begin{rema+}}{\end{rema+}}
\newenvironment{definition}{\begin{defi+}}{\end{defi+}}

\newcommand {\bC} {\mathbb {C}}

\newcommand {\bR} {\mathbb {R}}
\newcommand {\bZ} {\mathbb {Z}}

\newcommand {\la} {\lambda}
\newcommand {\calS} {\mathcal {S}}
\newcommand {\calD} {\mathcal {D}}
\newcommand {\calH} {\mathcal {H}}
\newcommand {\calC} {\mathcal {C}}
\newcommand {\calR} {\mathcal {R}}

\begin{document}

\numberwithin{equation}{section}

\title{Hyperbolic polynomials and spectral order}
\author[J.~Borcea]{Julius Borcea}
\address{Department of Mathematics, Stockholm University, SE-106 91 Stockholm,
   Sweden}
\email{julius@math.su.se}
\author[B.~Shapiro]{Boris Shapiro}
\address{Department of Mathematics, Stockholm University, SE-106 91 Stockholm,
   Sweden}
\email{shapiro@math.su.se}


\begin{abstract}
The spectral order on $\bR$ induces a partial ordering on the manifold 
$\calH_{n}$ of monic hyperbolic polynomials of degree $n$. We show that the 
semigroup $\tilde{\calS}$ generated by differential operators of the form 
$\left(1-\la \frac{d}{dx}\right)e^{\la \frac{d}{dx}}$, $\la \in \bR$, 
acts on the poset $\calH_{n}$ in an order-preserving fashion. We also show 
that polynomials 
in $\calH_{n}$ are global minima of their respective $\tilde{\calS}$-orbits 
and we conjecture 
that a similar result holds even for complex polynomials. Finally, we show 
that only those pencils of polynomials in $\calH_{n}$ which are of logarithmic 
derivative type satisfy a certain local minimum property for the 
spectral order.
\end{abstract}

\maketitle

\section*{Introduction and Main Results}

Given a complex polynomial $P$ of degree $n$ we define $Z(P)$ to be the 
unordered $n$-tuple consisting of the zeros of $P$, where each zero occurs as 
many times as its multiplicity. We denote by $\Re Z(P)$ the (unordered) 
$n$-tuple consisting of the real parts of the points in $Z(P)$. The 
polynomial $P$ is said to be {\em hyperbolic} if all its zeros are real. Note 
that in this case $\Re Z(P)=Z(P)$. A hyperbolic 
polynomial whose zeros are simple is called {\em strictly hyperbolic}. 

The main purpose of this paper is to study the behaviour of the $n$-tuple 
$Z(P)$ under the action of certain semigroups of differential operators. For 
this we shall use the following fundamental result from the theory of 
stochastic majorizations:

\begin{theorem}\label{spec}
Let $X=(x_{1},x_{2},\ldots,x_{n})^{t}$ and $Y=(y_{1},y_{2},\ldots,y_{n})^{t}$ 
be two unordered $n$-tuples of vectors in $\bR^{k}$. The following conditions 
are equivalent:
\begin{enumerate}
\item For any convex function $f: \bR^k \to \bR$ one
has $\sum_{i=1}^{n}f(x_{i})\le \sum_{i=1}^{n}f(y_{i})$.
\item There exists a doubly stochastic $n\times n$ matrix $A$ such that
$\tilde  X=A\tilde  Y$, where $\tilde X$ and $\tilde Y$ are 
$n\times k$ matrices
obtained by some (and then any) ordering of the vectors in $X$ and $Y$.
\end{enumerate}
\end{theorem}
If the conditions of Theorem~\ref{spec} are satisfied then we say that 
$X$ is \emph{majorized} by $Y$ or $X$ is {\it less than} $Y$ 
{\it in the spectral order}, and write $X\prec Y$.
Theorem~\ref{spec} is due to Schur and to Hardy, 
Littlewood and P\'olya in the one-dimensional case (\cite{HLP}), and to 
Sherman in the multivariate case (\cite{S}). These cases are also known 
as classical and multivariate majorization, respectively. One can easily 
check that $X\prec Y$ implies that 
$\sum_{i=1}^nx_{i}=\sum_{i=1}^ny_{i}$.

Let $\calH_{n}$ 
denote the manifold of monic hyperbolic polynomials of degree $n$. We may 
view $(\calH_{n},\preccurlyeq)$ as a partially ordered set, where the 
ordering relation 
$\preccurlyeq$ is induced by the spectral order on $n$-tuples of real numbers 
(cf.~Theorem \ref{spec}). Thus, if $P,Q\in \calH_{n}$ then 
$P\preccurlyeq Q$ if and 
only if $Z(P)\prec Z(Q)$. Note that although the spectral order is 
only a preordering on $n$-tuples of points in $\bR$, Birkhoff's theorem 
(\cite[Theorem 2.A.2]{MO}) implies that it actually 
induces a partial ordering on $\calH_{n}$.

Define the following semigroups of differential operators:
$$\calS=\Big\langle 1-\la \frac{d}{dx}\,\Big|\,\la \in \bR\Big\rangle,\quad 
\tilde{\calS}=\Big\langle \Big(1-\la \frac{d}{dx}\Big)e^{\la \frac{d}{dx}}
\,\Big|\,\la \in \bR\Big\rangle.$$
Note that $\tilde{\calS}$ is the (largest) subsemigroup of 
$\calS\times \langle e^{\mu \frac{d}{dx}}\mid \mu\in \bR\rangle$ which 
consists of operators that 
preserve the averages of the zeros of polynomials in $\calH_{n}$. 
The operator $\left(1-\la \frac{d}{dx}\right)e^{\la \frac{d}{dx}}$, 
$\la \in \bR$, will be denoted by $D_{\la}$ throughout this paper. It follows 
from the well-known Hermite-Poulain theorem (see \cite{O}) that the 
semi\-groups $\calS$ and $\tilde{\calS}$ act on $\calH_{n}$. Our first main 
result asserts that these semigroups act in fact on $\calH_{n}$ in an 
order-preserving fashion:

\begin{theorem}\label{semigp}
Let $P,Q\in \calH_{n}$ be such that $P\preccurlyeq Q$. Then $P+\la 
P'\preccurlyeq Q+\la Q'$ for any $\la \in \bR$.
\end{theorem}

We point out an interesting consequence of Theorem \ref{semigp}: 

\begin{corollary}
If $P,Q\in \calH_{n}$ are such that $P\preccurlyeq Q$ then 
$n^{-1}P'\preccurlyeq  n^{-1}Q'$.
\end{corollary}

The next theorem shows that any polynomial in $\calH_{n}$ is the global 
minimum of its orbit under the action of the semigroup $\tilde{\calS}$.

\begin{theorem}\label{orbit}
If $P\in \calH_{n}$ then $P\preccurlyeq D_{\la}P$ for any $\la \in \bR$.
\end{theorem}

A well-known theorem of Obreschkoff (see \cite{O}) states that if $P$ and $Q$ 
are real polynomials then the linear pencil of polynomials $P+\la Q$, 
$\la \in \bR$, consists 
of hyperbolic polynomials if and only if $P$ and $Q$ are hyperbolic and those 
of their zeros which are not common separate each other. The following 
converse to Theorem \ref{orbit} shows that real lines in $\calH_{n}$ 
of the form $P+\la P'$ are characterized by a 
local minimum property 
with respect to the partial ordering $\preccurlyeq$ on $\calH_{n}$.

\begin{theorem}\label{local}
Let $P\in \calH_{n}$, $\la \in \bR$, and let $Q$ be a complex polynomial of 
degree at most $n-1$. Set $R_{\la}(x)=P(x+\la)-\la Q(x+\la)$. If 
$R_{\la}\in \calH_{n}$ and $R_{0}\preccurlyeq R_{\la}$ for all small 
$\la \in \bR$ then $Q=P'$.
\end{theorem}

We also obtain a generalization of Theorem \ref{orbit} which shows that 
real lines in $\calH_{n}$ of the form $P+\la P'$ satisfy in fact a global 
monotony property:

\begin{theorem}\label{global}
If $P\in \calH_{n}$ and $\la_{1},\la_{2}\in \bR$ are such that 
$\la_{1}\la_{2}\ge 0$ and $|\la_{1}|\le |\la_{2}|$ then 
$D_{\la_{1}}P\preccurlyeq D_{\la_{2}}P$.
\end{theorem}

Finally, we show that real lines in $\calH_{n}$ 
of the form $P+\la P'$ satisfy an inequality {\it \`a la} G\aa rding 
(cf.~\cite{G}):

\begin{theorem}\label{convex}
If $P\in \calH_{n}$ then $\bR\ni \la \mapsto \max Z(D_{\la}P)$ is a convex 
function with a global minimum at $\la=0$ while $\bR\ni \la \mapsto 
\min Z(D_{\la}P)$ is a concave function with a global maximum at $\la=0$.
\end{theorem}

\begin{remark}
It follows from Theorem \ref{convex} that the so-called spread function 
$\bR\ni \la \mapsto \max Z(D_{\la}P)-\min Z(D_{\la}P)$ is a convex 
function with a global minimum at $\la=0$.
\end{remark}

The structure of the paper is as follows: in \S 1 we sketch the proofs of 
our main results and in \S 2 we present further questions and conjectures. 
The complete proofs will appear elsewhere.

\section{Outline of the proofs}

One of the key ingredients in the proofs of Theorems 
\ref{semigp}-\ref{convex} is the following 
criterion for classical majorization due to Hardy, Littlewood and P\'olya 
(cf.~\cite{HLP}). We should mention that there are no known analogues of this 
criterion for multivariate majorization.

\begin{theorem}\label{crit}
Let $X=(x_{1}\le x_{2} \le\ldots \le x_{n})\subset \bR$ and
$Y=(y_{1}\le y_{2} \le \ldots \le y_{n})$ be two $n$-tuples of 
real numbers. Then $X\prec Y$ if and only if the $x_{i}$'s and the 
$y_{i}$'s 
satisfy the following conditions:
   \begin{enumerate}
    \item $\sum_{i=1}^{n}x_{i}=\sum_{i=1}^{n}y_{i}$;
    \item $\sum_{i=0}^{k}x_{n-i}\le \sum_{i=0}^{k}y_{n-i}$ for $0\le k\le n-2$.
   \end{enumerate}
\end{theorem}

The proof of Theorem \ref{semigp} is based on several auxiliary results. 
Let us first make the following definition:

\begin{definition}
Let $P(x)=\prod_{i=1}^{n}(x-x_{i})\in \calH_{n}$, $n\ge 2$, and 
$1\le k<l\le n$. 
Assume that $x_{i}\le x_{i+1}$, $1\le i\le n-1$, and that $x_{k}\neq x_{l}$. 
Let further $t\in \,\,]0,\frac{x_{l}-x_{k}}{2}]$ and define $Q\in \calH_{n}$ 
to be 
the polynomial with zeros $y_{i}$, $1\le i\le n$, where $y_{k}=x_{k}+t$, 
$y_{l}=x_{l}-t$, and $y_{i}=x_{i}$, $i\neq k,l$. The polynomial $Q$ is 
called the {\em contraction of} $P$ {\em of type} $(k,l)$ 
{\em and coefficient} $t$. The contraction is called {\em simple} if $l=k+1$ 
and it is called {\em non-degenerate} if $t\neq \frac{x_{l}-x_{k}}{2}$.
\end{definition}

\begin{proposition}\label{contrac}
Let $P,Q\in \calH_{n}$ be two distinct strictly hyperbolic polynomials such 
that 
$P\preccurlyeq Q$. Then there exists a finite sequence 
$P_{1},\ldots,P_{m}\in \calH_{n}$ such that $P_{1}=Q$, $P_{m}=P$, and 
$P_{i+1}$ 
is a simple non-degenerate contraction of $P_{i}$, $1\le i\le m-1$.
\end{proposition}

\begin{remark}
Proposition \ref{contrac} is true even for polynomials with multiple zeros if 
the non-degeneracy condition is omitted.
\end{remark}

\begin{proposition}\label{simple}
Theorem \ref{semigp} is true if $P$ and $Q$ are strictly hyperbolic 
polynomials and $P$ is a simple (non-degenerate) contraction of $Q$.
\end{proposition}

From Propositions \ref{contrac} and \ref{simple} we deduce that Theorem 
\ref{semigp} is true in the generic case when $P$ and $Q$ have simple zeros. 
If this is not the case then we let $x_{i}$, $1\le i\le n$, and $y_{i}$, 
$1\le i\le n$, denote the zeros of $P$ and $Q$, respectively, and we choose 
an arbitrary positive number $\varepsilon$. Let $P_{\varepsilon}$ and 
$Q_{\varepsilon}$ be the polynomials with zeros $x_{i}-(n-i)\varepsilon$, 
$1\le i\le n-1$, $x_{n}+\frac{n(n-1)}{2}\varepsilon$, and 
$y_{i}-(n-i)\varepsilon$, $1\le i\le n-1$, 
$y_{n}+\frac{n(n-1)}{2}\varepsilon$, respectively. Note that 
$P_{\varepsilon}$ and $Q_{\varepsilon}$ are strictly hyperbolic and that 
$P_{\varepsilon}\preccurlyeq Q_{\varepsilon}$. The above arguments imply that 
$P_{\varepsilon}+\la P_{\varepsilon}'\preccurlyeq Q_{\varepsilon}+
\la Q_{\varepsilon}'$ for any $\la \in \bR$. Theorem \ref{semigp} now 
follows by letting $\varepsilon \rightarrow 0$.

Let $P(x)=\prod_{i=1}^{n}(x-x_{i})\in \calH_{n}$, $n\ge 2$, and $\la \in \bR$. 
Assume that $x_{i}<x_{i+1}$, $1\le i\le n-1$, and let $x_{i}(\la)$, 
$1\le i\le n$, denote the zeros of $D_{\la}P$. If these are labeled so that 
$x_{i}(0)=x_{i}$, $1\le i\le n$, then one can show that 
$x_{i}(\la)<x_{i+1}(\la)$ and that by varying $x_{n}$ and keeping $\la$ fixed 
each $x_{i}(\la)$, $1\le i\le n-1$, is an increasing function of $x_{n}$. This 
makes it possible to prove Theorem \ref{orbit} by induction on $n$ in 
the generic case. If $P\in \calH_{n}$ has multiple zeros we notice that 
$\left(1-\varepsilon \frac{d}{dx}\right)^{n-1}\!P$ has simple zeros for any 
$\varepsilon\neq 0$. Since Theorem \ref{orbit} is true for the latter 
polynomials, we get that Theorem \ref{orbit} holds in the general case  
by letting $\varepsilon \rightarrow 0$.

To prove Theorem \ref{local} we first use Theorem \ref{crit} in order to 
show that if $P$ is strictly 
hyperbolic and $Q$ satisfies the assumptions of Theorem \ref{local} then 
$Q(x_{i})=P'(x_{i})$, where $x_{i}$, $1\le i\le n$, are the zeros of $P$. 
Using again the fact that $\left(1-\varepsilon \frac{d}{dx}\right)^{n-1}\!P$ 
has simple zeros if $P\in \calH_{n}$ and $\varepsilon\neq 0$ and also that 
the operator $\left(1-\varepsilon \frac{d}{dx}\right)^{n-1}$ preserves the 
ordering 
on $\calH_{n}$ (cf.~Theorem \ref{semigp}) we deduce that Theorem \ref{local} 
holds for any $P\in \calH_{n}$.

Let now $P(x)=\prod_{i=1}^{n}(x-x_{i})\in \calH_{n}$, $n\ge 2$, 
$\la \in \bR$,  
and assume that $x_{i}<x_{i+1}$, $1\le i\le n-1$. Denote the zeros of 
$D_{\la}P$ by $x_{i}(\la)$, $1\le i\le n$, those of $P'$ by $w_{j}$, 
$1\le j\le n-1$, and those of $D_{\la}P'$ by $w_{j}(\la)$, $1\le j\le n-1$. 
If these are labeled so that $x_{i}(0)=x_{i}$, $1\le i\le n$, and 
$w_{j}(0)=w_{j}$, $1\le j\le n-1$, then one can show that 
$x_{i}(\la)<w_{i}(\la)<x_{i+1}(\la)$, $1\le i\le n-1$, and that for any 
$\la>0$ and $1\le m\le n-1$ one has
\begin{equation*}
\sum_{i=1}^{m}x_{i}'(\la)=\la \sum_{j=1}^{n-1}\frac{P''(w_{j}(\la)+\la)}
{D_{\la}P''(w_{j}(\la))}\left(\sum_{i=1}^{m}\frac{1}{x_{i}(\la)-w_{j}(\la)}
\right)<0.
\end{equation*}
Thus $\la\mapsto \sum_{i=1}^{m}x_{i}(\la)$, $1\le m\le n-1$, are decreasing 
functions on $]0,\infty[$, which combined with Theorem \ref{crit} proves 
Theorem \ref{global} for strictly hyperbolic polynomials (the case $\la<0$ is 
similar). The same density arguments as those used for Theorem \ref{local} 
show that Theorem \ref{global} is true for all $P\in \calH_{n}$.

The second part of the statement in Theorem \ref{convex} follows from 
the first part by noticing that $P(x)\in \calH_{n}$ if and 
only if $(-1)^nP(-x)\in \calH_{n}$. Keeping the same notations as above, 
one can show that if $P\in \calH_{n}$ is strictly hyperbolic then
\begin{equation*}
x_{n}''(\la)=2(x_{n}'(\la)+1)^2\sum_{i=1}^{n-1}\frac{w_{i}-x_{i}(\la)-\la}
{(x_{n}(\la)+\la-w_{i})(x_{n}(\la)-x_{i}(\la))}>0
\end{equation*}
for any $\la \in \bR$, which proves Theorem \ref{convex} in the generic case. 
If $P$ has multiple zeros then one can use strictly hyperbolic polynomials 
of the form $\left(1-\frac{1}{s}\frac{d}{dx}\right)^{n-1}\!P$, $s\in \bZ_{+}$,
 in order to approximate the function $\la \mapsto \max Z(D_{\la}P)$ 
uniformly on compact intervals by convex $C^2$-functions. This proves 
Theorem \ref{convex}. 

\section{Remarks and open questions}

The manifold $\calC_{n}$ of monic complex polynomials of degree $n$ is a 
natural 
context for discussing possible extensions of the above results to the 
complex case. By analogy with the hyperbolic case we may view 
$(\calC_{n},\preccurlyeq)$ as a partially ordered set, where the ordering 
relation $\preccurlyeq$ is now induced by the spectral order on  
$n$-tuples of 
vectors in $\bR^2$ (cf.~Theorem \ref{spec}). This means that the zero sets 
of polynomials in $\calC_{n}$ are viewed as subsets of $\bR^{2}$ and that if 
$P,Q\in \calC_{n}$ then $P\preccurlyeq Q$ if and only if $Z(P)\prec 
Z(Q)$.

The following example shows that if the partial ordering $\preccurlyeq$ on 
$\calC_{n}$ is 
defined as above then one cannot expect a complex analogue of 
Theorem \ref{orbit}.

\begin{proposition}\label{counter1}
Let $P(z)=z^n-1$ and $\la \in \bC$. If $n\ge 3$ and $|\la|$ is small enough 
then $D_{\la}P$ and $P$ are 
incomparable with respect to the partial ordering $\preccurlyeq$ on 
$\calC_{n}$.
\end{proposition}

We also note that the results of the previous section (for the 
hyperbolic case) are valid only for real values of the parameter $\la$:

\begin{proposition}\label{counter2}
Let $P\in \calH_{n}$ be a strictly hyperbolic polynomial and let 
$\la \in \bC\setminus \bR$. If $n\ge 2$ and $|\la|$ is sufficiently small 
then $D_{\la}P$ and $P$ are incomparable with respect to the partial 
ordering $\preccurlyeq$ on $\calC_{n}$.
\end{proposition}

These examples suggest that complex 
generalizations of Theorem \ref{orbit} -- if any -- should involve only 
classical majorization and real values of the parameter $\la$. Based on 
extensive numeric calculations, we make the following

\begin{conjecture}\label{conj1}
If $P\in \calC_{n}$ then $\Re Z(P)\preccurlyeq \Re Z(D_{\la}P)$ for any $\la 
\in \bR$. 
\end{conjecture}

We end with a few questions related to the semigroups $\calS$ and 
$\tilde{\calS}$. Let $\calR_{n}$ denote the set of all monic real polynomials 
of degree $n$. There is reason to believe that any two $\calS$-orbits in 
$\calR_{n}$ have a non-empty intersection:

\begin{conjecture}\label{conj2}
If $P_{1},P_{2}\in \calR_{n}$ then there exist differential operators 
$\Lambda_{1},\Lambda_{2}\in \calS$ such that $\Lambda_{1}P_{1}
=\Lambda_{2}P_{2}$.
\end{conjecture}

If true, Conjecture \ref{conj2} would imply in particular that 
$\calS P\,\cap\,\calH_{n}\neq 
\emptyset$ for any $P\in \calR_{n}$, which would answer in the affirmative a 
question of I.~Krasikov.

Let finally $P\in \calH_{n}$ and set $P_{\preccurlyeq}=\{Q\in \calH_{n}\mid 
P\preccurlyeq 
Q\}$. One can easily check that if $P$ is strictly hyperbolic and $n\ge 3$ 
then $\tilde{\calS}P\subsetneq P_{\preccurlyeq}$. It would be interesting to 
know 
whether there exists a (semi)group $\calD$ of differential operators (not 
necessarily with constant coefficients) such that 
$\calD\supsetneq \tilde{\calS}$ and $P_{\preccurlyeq}=\calD P$ for any 
$P\in \calH_{n}$. This would give a 
completely new way of describing classical majorization.

\medskip

\noindent
{\bf Acknowledgements.} The results announced in this paper and their proofs 
were all obtained by the first author (JB). The second author (BS) only made 
some computer tests in connection with Remark 1 and Theorem~\ref{orbit}. Both 
authors are grateful to Harold Shapiro and Ilia Krasikov for stimulating 
discussions.

\end{document}